\documentclass[12pt]{article}

\usepackage{amsmath,amsthm,amssymb,epsfig,xargs,hyperref}
\usepackage{colortbl}
\providecommand{\mk}{\cellcolor[gray]{.6}}

\usepackage[margin=1.5cm]{geometry}

\newtheorem{theorem}{Theorem}

\usepackage[colorinlistoftodos,prependcaption,textsize=tiny]{todonotes}
\newcommandx{\unsure}[2][1=]{\todo[linecolor=red,backgroundcolor=red!25,bordercolor=red,#1]{#2}}
\newcommandx{\change}[2][1=]{\todo[linecolor=blue,backgroundcolor=blue!25,bordercolor=blue,#1]{#2}}
\newcommandx{\info}[2][1=]{\todo[linecolor=OliveGreen,backgroundcolor=OliveGreen!25,bordercolor=OliveGreen,#1]{#2}}
\newcommandx{\improvement}[2][1=]{\todo[linecolor=Plum,backgroundcolor=Plum!25,bordercolor=Plum,#1]{#2}}
\newcommandx{\thiswillnotshow}[2][1=]{\todo[disable,#1]{#2}}

\def\dfrac#1#2{\lower0.15ex\hbox{\large$\frac{#1}{#2}$}} % bdm's tricky fracs

\renewcommand{\leq}{\leqslant}

\renewcommand{\le}{\leqslant}

\def\eref#1{$(\ref{#1})$}

\def\sref#1{\S$\ref{#1}$}

\def\cref#1{Corollary~$\ref{#1}$}

\def\T{{\ensuremath{\cal T}}}
\def\F{{\ensuremath{\cal F}}}
\def\U{{\ensuremath{\cal U}}}
\def\I{{\ensuremath{\cal I}}}

\def\sym{\mathcal{S}}

\title{
Perfect 1-factorisations of $K_{16}$
}

\author{
  Michael J. Gill\thanks{Research supported by an Australian Government Research Training Program (RTP) Scholarship.}, \ Ian M.\ Wanless\thanks{%
Research supported by ARC grant DP150100506.}\\
\small School of Mathematics\\[-0.75ex]
\small Monash University\\[-0.75ex]
\small Vic 3800, Australia\\[-0.75ex]
\small\tt \{michael.gill,ian.wanless\} \ @monash.edu
}

\date{}

\begin{document}

\maketitle

\begin{abstract}
We report the results of a computer enumeration that found that there
are 3155 perfect 1-factorisations (P1Fs) of the complete graph $K_{16}$.
Of these, 89 have a non-trivial automorphism group (correcting an earlier
claim of 88 by Meszka and Rosa \cite{MR03}).

We also (i) describe a new invariant which distinguishes between the
P1Fs of $K_{16}$, (ii) observe that the new P1Fs produce no atomic
Latin squares of order 15 and (iii) record P1Fs for a number of large
orders that exceed prime powers by one.
\end{abstract}

\section{Introduction}\label{s:intro} 

A $k$-factor in a graph is a $k$-regular spanning subgraph. In
particular, a $1$-factor (also called a perfect matching) is a set of
edges in a graph which cover every vertex exactly once. A
$1$-factorisation is a collection of 1-factors which partitions the
edges of the graph.  Suppose that we have a 1-factorisation $\F$ of
some graph. The union of any two distinct 1-factors in $\F$ is a
2-factor, that is, a collection of cycles. If, regardless of which two
1-factors we choose in $\F$, their union is a single cycle, then we say
that $\F$ is perfect.  Throughout, we will use the abbreviation P1F
for \emph{perfect $1$-factorisation}. For background reading on P1Fs,
including definitions of terms not included here, we refer to
\cite{Sea91,Wal97}.

The P1Fs of complete graphs up to $K_{14}$ have been known for some
time \cite{DG96}. The main purpose of this note is to report on a
computer enumeration of the next case, namely P1Fs of $K_{16}$.  We
also check the Latin squares associated with the new P1Fs
(\sref{s:ALS}) and find they are not atomic, discuss invariants which
distinguish non-isomorphic P1Fs (\sref{s:invar}), and record a number
of new P1Fs for large orders that are one more than a prime power
(\sref{s:largeord}).

Kotzig's P1F conjecture famously asserts that P1Fs of $K_n$ exist for all even $n$.
Only three infinite families are known \cite{BMW06}. They cover all orders of the 
form $n=p+1$ or $n=2p$ where $p$ is an odd prime. Sporadic constructions including
those reported in \sref{s:largeord}, together with \cite{Pik19,Sea91,cyclatom,Wol09}, demonstrate
existence for $n\le56$ and the following orders:
\begin{align*}
&126,170,244,344,530,730,1332,1370,1850,2198,2810,3126,4490,6860,
6890,11450,11882,12168,\\
&15626,16808,22202,24390,24650,26570,29792,29930,32042,38810,
44522,50654,51530,52442,63002,\\
&72362,76730,78126,79508,103824,
148878,161052,205380,226982,300764,357912,371294,493040,\\
&571788,1030302,1092728,1225044,1295030,2048384,2248092,2476100,2685620,3307950,
3442952,\\
&4330748,4657464,5735340,6436344,6967872,7880600,9393932,11089568,
11697084,13651920,\\
&15813252,18191448,19902512,22665188.
\end{align*}
Kotzig's P1F conjecture remains open for other orders. The smallest unsolved cases are\\
$\{64, 66, 70, 76, 78, 88, 92, 96, 100\}$.

\section{The catalogue}\label{s:computation}

We labelled the 16 vertices of our complete graph {\tt a,b,$\dots$,p},
and assumed without loss of generality that every P1F contained the
two factors
\begin{equation}\label{e:F1F2}
  \hbox{{\tt $F_1=\{$ab,cd,ef,gh,ij,kl,mn,op$\}$} and {\tt $F_2=\{$ap,bc,de,fg,hi,jk,lm,no$\}$}}.
\end{equation}
Next, we generated all 1-factors of
$K_{16}$ containing the edge {\tt ac}. Up to isomorphism, this gave us
1647 options for the initial three 1-factors $F_1,F_2,F_3$.  In each
isomorphism class we chose $F_3$ to be the lexicographically least
option.  Then for each of these isomorphism representatives we did the
following.  First, we generated and stored the set $\T$ of all
1-factors that were compatible with $F_1,F_2,F_3$, in the sense that
$F_i\cup F$ is a Hamilton cycle for all $F \in \T$ and
$i\in\{1,2,3\}$. The number of 1-factors we needed to store at this
stage depended on the choice of $F_1, F_2, F_3$ and ranged from
56\,816 to 59\,312. We also computed and stored which pairs of these
1-factors in $\T$ were compatible (made a Hamilton cycle).  A list of
``active'' 1-factors was maintained, which initially consisted of $\T$.

We then used backtracking to add one 1-factor at a time from the
active list.  Each new factor was forced to include a particular edge
that was judged to be the best option of the edges remaining to be
used. The criterion was to choose the edge that was contained in the
fewest active 1-factors, thereby limiting the branching that our
search undertook. Having chosen a new 1-factor, we checked that no
relabelling of the current partial 1-factorisation would give a
lexicographically smaller set of 3 initial factors than the one
currently being processed. Assuming that to be the case, we updated
the list of active 1-factors by removing any which were incompatible with
the new 1-factor, and continued the search. Any P1Fs that
we found were canonically labelled (see \sref{s:invar}), then output.
The result was

\begin{theorem}\label{t:3155}
There are $3155$ P1Fs of $K_{16}$ up to isomorphism. Of these, $89$ have
a non-trivial automorphism group.
\end{theorem}

The 3155 P1Fs can be downloaded from \cite{catalog}. Note that the same
catalogue has been independently and concurrently generated by Mariusz Meszka \cite{Mes20}.
The P1Fs with non-trivial automorphism groups all had cyclic automorphism
groups, with generators as follows:
\begin{itemize}
\item A unique starter induced P1F, with automorphism group generated by
a permutation with cycle type $15^11^1$.

\item A unique even-starter induced P1F, with automorphism group generated by
a permutation with cycle type $14^11^2$.

\item 4 P1Fs with automorphism group generated by
a permutation with cycle type $7^21^2$.

\item 5 P1Fs with automorphism group generated by
a permutation with cycle type $5^31^1$.

\item 19 P1Fs with automorphism group generated by
a permutation with cycle type $3^51^1$.

\item 59 P1Fs with automorphism group generated by
a permutation with cycle type $2^71^2$.

\end{itemize}

These numbers agree with an earlier enumeration by Meszka and Rosa \cite{MR03},
{\em except} that we found one extra P1F with automorphism group of order 7.
That P1F can be obtained by developing the following two 1-factors under the
permutation {\tt(abcdefg)(hijklmn)}:
\begin{equation}\label{e:newaut2}
\hbox{\tt $\{$ab, cg, do, em, fi, hp, jl, kn$\}$,
  $\{$ac, bk, dj, ei, fp, gl, ho, mn$\}$,}
\end{equation}
together with the 1-factor {\tt$\{$ah, bi, cj, dk, el, fm, gn, op$\}$}.

\section{Connections with Latin squares}\label{s:ALS}

In this section we explore whether the newly discovered P1Fs of $K_{16}$
produce interesting Latin squares. First we need to give some definitions.

A {\em Latin square} of order $n$ is an $n\times n$ matrix in which
each row and column is a permutation of some (fixed) symbol set of
size $n$, say $\{1,2,\dots,n\}$.  It is often useful to think of a
Latin square of order $n$ as a set of $n^2$ triples $(i,j,k)$ where
$k$ is the symbol that appears in cell $(i,j)$.  Each Latin square then has
six {\em conjugate} squares obtained by (uniformly) permuting the
coordinates of each triple.  Writing Latin squares of order $n$ as
sets of triples also provides for a natural action of $\sym_n\wr\sym_3$ on
the Latin squares.  The \emph{species} (sometimes called \emph{main
  class}) of a Latin square is its orbit under this action.

A \emph{Latin subrectangle} is a rectangular submatrix of a Latin
square $L$ in which the same symbols
occur in each row.  If $R$ is a $2\times\ell$ Latin subrectangle of
$L$, and $R$ is minimal in that it contains no $2\times\ell'$ Latin
subrectangle for $2\leq\ell'<\ell$, then we say that $R$ is a \emph{row
cycle of length $\ell$}.  
Column cycles and symbol cycles can be defined similarly, and the
operations of conjugacy on $L$ interchange these objects.

A Latin square of order $n$ is {\it row-Hamiltonian} if every row
cycle has length $n$, \emph{symbol-Hamiltonian} if every symbol cycle has
length $n$, and \emph{column-Hamiltonian} if every column cycle has length
$n$.  These three types of square are related by conjugacy.
A Latin square is {\it atomic} if it is row-Hamiltonian,
symbol-Hamiltonian and column-Hamiltonian.  In other words, a square
is atomic if all its conjugates are row-Hamiltonian. 

A $1$-factorisation $\F$ of $K_n$ is equivalent to a symmetric unipotent
Latin square $\U(\F)$, where we define $\U(\F)[i,j]=k$ if the edge $ij$ appears
in the $k$-th factor of $\F$, and $\U(\F)[i,i]=n$ for all $i$.
From $\U(\F)$ we can form $n$ idempotent symmetric Latin squares
of order $n-1$, by ``folding'' a row and the corresponding column onto
the main diagonal. To be precise, we form $\I(\F,j)$ from 
$\U(\F)$ by replacing the entry in cell $(i,i)$ by $\U(\F)[i,j]$ for
each $i$, and then deleting row $j$ and column $j$.
The 1-factorisation $\F$ is perfect if and only if for each $j$ the
corresponding Latin square $\I(\F,j)$ is symbol-Hamiltonian
\cite{syminherit}. There is some interest in whether $\I(\F,j)$
has the stronger property of being atomic.  Atomic Latin squares are
known to exist for some composite orders, but no example is known of
an atomic Latin square whose order is not a prime power
\cite{cyclatom}. This makes order 15, the smallest case where
existence of an atomic Latin square is open, particularly interesting.

Each P1F $\F$ corresponds to $k$ species of Latin squares $\I(\F,j)$,
where $k$ is the number of orbits of the automorphism group of $\F$. Hence
P1Fs with automorphism groups generated by permutations of cycle type
$15^11^1,14^11^2,7^21^2,5^31^1,3^51^1,2^71^2,1^{16}$, respectively produce
$2,3,4,4,6,9,16$ species of Latin squares. Hence, overall our 3155 P1Fs of
$K_{16}$ produce 49\,742 species containing symmetric symbol-Hamiltonian
Latin squares. Checking a representative of each species, we established:

\begin{theorem}
There exists no symmetric atomic Latin square of order $15$.
\end{theorem} 

The nearest that we got to an atomic Latin square was $N_{15}$, where
\[
N_{15}=\left[
\begin{array}{ccccccccccccccc}
 1&12& 2& 7&10&\mk 8& 5&15&\mk13&\mk 3& 6& 4&14&11& 9\\
12& 2&13& 3& 8&11& 9& 6&15&14& 4& 7& 5& 1&10\\
 2&13& 3&14& 4& 9&12&10& 7&15& 1& 5& 8& 6&11\\
 7& 3&14& 4& 1& 5&10&13&11& 8&15& 2& 6& 9&12\\
10& 8& 4& 1& 5& 2& 6&11&14&12& 9&15& 3& 7&13\\
 8&11& 9& 5& 2& 6& 3& 7&12& 1&13&10&15& 4&14\\
 5& 9&12&10& 6&\mk 3& 7& 4&\mk 8&\mk13& 2&14&11&15& 1\\
15& 6&10&13&11& 7& 4& 8& 5& 9&14& 3& 1&12& 2\\
13&15& 7&11&14&12& 8& 5& 9& 6&10& 1& 4& 2& 3\\
 3&14&15& 8&12& 1&13& 9& 6&10& 7&11& 2& 5& 4\\
 6& 4& 1&15& 9&13& 2&14&10& 7&11& 8&12& 3& 5\\
 4& 7& 5& 2&15&10&14& 3& 1&11& 8&12& 9&13& 6\\
14& 5& 8& 6& 3&15&11& 1& 4& 2&12& 9&13&10& 7\\
11& 1& 6& 9& 7& 4&15&12& 2& 5& 3&13&10&14& 8\\
 9&10&11&12&13&14& 1& 2& 3& 4& 5& 6& 7& 8&15\\
% Or can use alphabetical format
%a&l&b&g&j&h&e&o&m&c&f&d&n&k&i\\
%l&b&m&c&h&k&i&f&o&n&d&g&e&a&j\\
%b&m&c&n&d&i&l&j&g&o&a&e&h&f&k\\
%g&c&n&d&a&e&j&m&k&h&o&b&f&i&l\\
%j&h&d&a&e&b&f&k&n&l&i&o&c&g&m\\
%h&k&i&e&b&f&c&g&l&a&m&j&o&d&n\\
%e&i&l&j&f&c&g&d&h&m&b&n&k&o&a\\
%o&f&j&m&k&g&d&h&e&i&n&c&a&l&b\\
%m&o&g&k&n&l&h&e&i&f&j&a&d&b&c\\
%c&n&o&h&l&a&m&i&f&j&g&k&b&e&d\\
%f&d&a&o&i&m&b&n&j&g&k&h&l&c&e\\
%d&g&e&b&o&j&n&c&a&k&h&l&i&m&f\\
%n&e&h&f&c&o&k&a&d&b&l&i&m&j&g\\
%k&a&f&i&g&d&o&l&b&e&c&m&j&n&h\\
%i&j&k&l&m&n&a&b&c&d&e&f&g&h&o\\
\end{array}\right].
\]
This square is derived from the unique even-starter induced P1F of
$K_{16}$.  It is one of the two species of Latin squares derived from
that P1F that inherit an automorphism group of order 14. Indeed, it is
clear that $N_{15}$ has an automorphism applying the cycle
$(1,2,\dots,14)$ simultaneously to rows, columns and symbols.
It is also clear that $N_{15}$ is not atomic, because it has a row-cycle of
length 3 (highlighted).  However, of the $\binom{15}{2}=105$ pairs of
rows in $N_{15}$, there are 91 that produce a Hamilton row-cycle.
The only pairs which fail are the 14 images under the automorphism group 
of the pair of rows containing the highlighted row-cycle.

\section{Invariants}\label{s:invar}

The literature contains several invariants that can be used for
distinguishing non-isomorphic P1Fs. The best known of these are the
train and tri-colour vector \cite{Wal97}. In our discussion we will
mention several \emph{complete invariants} for P1Fs of $K_{16}$. That is, invariants
which coincide for two P1Fs if and only if the P1Fs are isomorphic.

The train itself is a complete invariant for the P1Fs of $K_{16}$. However,
the indegree sequence of the train is not; it partitions the 3155
P1Fs into 3104 equivalence classes. The P1Fs that do not have
unique indegree sequences fall into 47 pairs and 2 triples. We now explicitly
give the two triples. The format is similar to that used in \eref{e:newaut2}, except
that we compress factors by omitting spaces and punctuation, and use only a space to
separate factors.
The first triple shares the indegree sequence $[573,784,336,86,19,2]$:

{\tt\noindent
  abcdefghijklmnop acbedgfihljmkonp adblcheofngmikjp aebnckdjfoglhpim afbicpdmeghojlkn\\
  agbdcnemfjhkiolp ahbpcodkelfmgijn aibmcedpflgkhnjo ajbgcidnehfklomp akbfcmdoepgnhjil\\
  albhcgdfejinkpmo ambjcfdhekgoipln anbocjdleifhgpkm aobkcldienfpgjhm apbcdefghijklmno

  \smallskip

  \noindent
  abcdefghijklmnop acbedgfihkjnlpmo adbkcienfhgpjmlo aebncgdfhlikjomp afbhcjdmeoglipkn\\
  agbfckdoejhpimln ahblcfdkemgoinjp aibmcpdnehfogkjl ajbdclepfmgihnko akbjcndielfpgmho\\
  albgcodjeifkhmnp ambpchdlekfngjio anbocmdhegfjilkp aobicedpflgnhjkm apbcdefghijklmno

  \smallskip

  \noindent
  abcdefghijklmnop acbedgfihljnkpmo adbfchejgmioknlp aebjcpdlfmgohnik afblckdoeghjinmp\\
  agbncmdkelfjhoip ahbpcndfeigjkmlo aibmcjdpenfhglko ajbkcgdiemfohpln akbocedhfngpiljm\\
  albicfdmeognhkjp ambdclehfkgijonp anbhcodjepflgkim aobgcidnekfphmjl apbcdefghijklmno

}

\medskip

Meanwhile, the second triple shares the indegree sequence $[584,765,338,94,18,1]$:

{\tt\noindent
abcdefghijklmnop acbedgfihmjokpln adbicpeoflgnhjkm aebgcldfhpiojmkn afbhckdnejgoimlp\\
agbjcodkeifmhlnp ahbpcgdienfkjlmo aiblcndoepfjgmhk ajbnchdmegfpilko akbfcjdpemgihnlo\\
albocmdjehfngkip ambkcedhfoglinjp anbdcfelgjhoikmp aobmcidlekfhgpjn apbcdefghijklmno

\smallskip

\noindent
abcdefghijklmnop acbedgfihljnkpmo adbmcjehfpgkioln aebhckdofngijlmp afbjcndleogmhpik\\
agbfcedmhkiljonp ahbncgdfeijpkmlo aibdcmejfhgoknlp ajbocidnepfkglhm akbpcodjegflhnim\\
albgcpdkemfohjin amblchdienfjgpko anbicldpekfmgjho aobkcfdhelgnipjm apbcdefghijklmno

\smallskip

\noindent
abcdefghijklmnop acbedgfihkjnlpmo adblciemfhgojpkn aebjcfdngmhoilkp afbdcoeigkhpjmln\\
agbpcndkehflimjo ahbmcjdpenfkgilo aibkcldheofjgnmp ajbgcedmfphlinko akbhcmdjelfngpio\\
albockdiepfmgjhn ambfchdoejgliknp anbicpdlegfohjkm aobncgdfekhmipjl apbcdefghijklmno 
}

\medskip

Among the pairs of P1Fs whose train indegree sequences coincide, there is (only) one that
involves two P1Fs whose automorphism groups have different order.
Specifically the following rigid P1F has indegree sequence $[598,748,332,102,18,2]$,
which coincides with one of the P1Fs with automorphism group of order 2 reported
in \cite{MR03}:

{\tt\noindent
% The one with grp of order 2
%abcdefghijklmnop acbedgfihkjmlonp adbfcjengohlikmp aebkcgdofnhjimlp afbicndpelgkhmjo\\
%agbochdmepfjilkn ahbmcidkegfojpln aibhcldfejgnkpmo ajblcedifmgphnko akbdcmehfpgliojn\\
%albgcpdjeofhinkm ambpcfdnekgihojl anbjcodleifkgmhp aobnckdhemflgjip apbcdefghijklmno
%\smallskip
%\noindent
abcdefghijklmnop acbedgfjhlinkomp adbmcpeiflgohkjn aebdclfogihmjpkn afbjcndkelgmhpio\\
agblcjdhekfnipmo ahbpcodnejfkglim aibhcfdoemgnjlkp ajbgcedmfhiklonp akbicgdpeofmhjln\\
albocmdjepfigkhn ambkchdfengpiljo anbfckdieghojmlp aobncidlehfpgjkm apbcdefghijklmno
}

\medskip

The train is a complete invariant for P1Fs of $K_{16}$, but comparing
trains is an instance of digraph isomorphism, which is not known to be
possible in polynomial time. That makes it desirable to identify
features of the train that are easily computed (like the indegree
sequence), but still have enough information to be a complete
invariant (unlike the indegree sequence). One such property is the
following. For each vertex $v$ of the train define $p(v)$ to be the
length of the shortest directed path from $v$ to any vertex $w$ that
is in a directed cycle. Since every vertex has outdegree 1, it is a
simple matter to construct a path starting at $v$ and following the
unique outgoing arc until we reach a vertex $w$ that we have
previously seen on the path. The distance from $v$ to the first
occurrence of $w$ is, by definition, $p(v)$. This also shows that
$p(v)$ is well defined. Note that if $v$ itself is in a cycle then
$p(v)=0$. Also, for a general $1$-factorisation the cycle involving
$w$ may be a loop. However, for a P1F of $K_n$ (where $n>2$) this
cannot happen, because loops only occur on vertices that are otherwise
isolated. To see this, consider two vertices $u=(\{a,b\},F)$ and
$v=(\{c,d\},G)$ of the train where $F,G$ are 1-factors and $a,b,c,d$
are vertices of the graph being factorised. If there was an arc $uv$
as well as a loop on $v$ that would mean (up to
swapping $c$ and $d$) that $ac,bd\in F$ and
$ab,cd\in G$, which produces a 4-cycle in $F\cup G$.

We found that for $i=0,1,2,3,4,5$ the counts of how many vertices have $p(v)=i$ together
form a complete invariant for the P1Fs of $K_{16}$. Counting only up to $i=4$ does not suffice,
since there is a single pair of (rigid) P1Fs whose first difference occurs at $i=5$. They are 

{\tt\noindent
abcdefghijklmnop acbedgfihmjokpln adbmcfeigjhlkonp aebdcnflgkhjiomp afbicjdkenglhpmo\\
agblcodhejfnipkm ahbkcpdjeofmgnil aibncldoepfjgmhk ajbfckdnemgiholp akbjcedlfpgohnim\\
albpchdiegfojmkn ambhcgdfekinjplo anbocidmehfkgpjl aobgcmdpelfhikjn apbcdefghijklmno\\
}
which has counts $[139,19,15,14,17,17]$, and 

{\tt\noindent
abcdefghijklmnop acbedgfihljokmnp adbicefmgjhoknlp aebkcndhfogpiljm afbpcjdlekgihnmo\\
agbjchdpemfkinlo ahbockdmeiflgnjp aibgcmdoepfjhkln ajbfcldnehgoikmp akbmcidfeoglhpjn\\
albhcpdiejfngmko amblcgdjenfhiokp anbdcoelfpgkhjim aobncfdkeghmipjl apbcdefghijklmno\\
}
which has counts $[139,19,15,14,17,22]$.

\medskip

Another well studied invariant of P1Fs is the tricolour vector \cite{Wal97}. It was already
reported in \cite{MR03} that the tricolour vector does not fully distinguish the
P1Fs of $K_{16}$ with non-trivial automorphism group. We found that it
partitions the set of all P1Fs of $K_{16}$ into 2320 equivalence classes.

Next we propose a new invariant that performs better than either the
train indegree sequence or tri-colour vector on the P1Fs of
$K_{16}$. It is based on the vertex cycles which were used for
switchings among 1-factorisations in \cite{KMOW14}.  These cycles
correspond exactly to the row-cycles in $\U(\F)$ other than the
row-cycles of length 2 that include entries on the main
diagonal. Finding the length of all row-cycles in $\U(\F)$ is easily
done in cubic time. The vector of tallies of how many row-cycles there
are of each length is a complete invariant for P1Fs of $K_{16}$.

Another invariant (that is presumably more discriminating in general),
is to count for each row of $\U(\F)$ how many row-cycles of each
length involve that row. The resulting list of $n$ vectors should then
be sorted lexicographically to accommodate relabelling of the vertices
of $\F$. We found that for the P1Fs of $K_{16}$ this invariant was
complete even if we only counted cycles of length 3 and 4. Counting
only the cycles of length 3 was not a complete invariant, but did
separate the 3155 P1Fs into 3102 equivalence classes.

The invariant involving $p(v)$ as described above can be calculated in polynomial
time, but not obviously in cubic time. By comparison,
the indegree sequence of the train, tri-colour vector and our new
invariant based on vertex cycle lengths are all easily computed in
time cubic in the order of the graph being factorised. 
The indegree sequence and tri-colour vector
are not complete invariants for P1Fs of $K_{16}$ and
it seems likely that the vertex cycle lengths and path lengths will not be a complete
invariant when $n$ gets larger.  However, it is possible in polynomial
time to canonically label a P1F, which does provide a complete
invariant for every order. We may define our canonical form to contain
the factors $F_1$ and $F_2$ from \eref{e:F1F2}, and to be the
lexicographically least possibility under that assumption.  Starting
with any P1F, there are only quadratically many choices for the
factors that will become $F_1$ and $F_2$. For each of those choices,
there are only linearly many ways that the two factors can map to $F_1$
and $F_2$ (given that their union is a single cycle). Thus we can check
all relabellings in polynomial time and choose the canonical one.
Once we have a canonical relabelling, isomorphism checking is a triviality.
The fact that isomorphism testing for P1Fs can by done in polynomial
time has long been known \cite{CC80,MR03}. The P1Fs given explicitly
above, except the one given in \eref{e:newaut2}, have been given in the
canonical form just described. All P1Fs in the catalogue \cite{catalog} are
also in this form.

\section{New orders of P1Fs}\label{s:largeord}

The paper \cite{cyclatom} gave constructions for many new orders of P1Fs
using the quotient coset starter method. Subsequent to publishing that paper
the author discovered many more examples and put them on his webpage, where
they have been visible for more than 12 years. In the interests of committing 
them to a more permanent part of the literature, we take this opportunity
to record them here. The P1Fs of order 1092728 and 1225044 were apparently
previously found by Volker Leck, but we are not aware of them having been published.

The following produce P1Fs for complete graphs of orders $q+1=p^3+1$,
where $p$ is prime. They can be considered as extra entries for Table
6 in \cite{cyclatom}, and that paper should be consulted for the
meaning of the notation.

$p=101$, $q=1030301$, $\zeta(x)=x^3+x+3$, $\tilde{c}=[813092,759910,233271,3]$

$p=103$, $q=1092727$, $\zeta(x)=x^3+x+4$, $\tilde{c}=[828376,896]$

$p=107$, $q=1225043$, $\zeta(x)=x^3+x+9$, $\tilde{c}=[1107573,151]$

$p=109$, $q=1295029$, $\zeta(x)=x^3+x+6$, $\tilde{c}=[271574,645911,1082655,4]$

$p=127$, $q=2048383$,  $\zeta(x)=x^3+x+15$, $\tilde{c}=[840749,23]$

$p=131$, $q=2248091$,  $\zeta(x)=x^3+x+3$, $\tilde{c}=[2096100,298]$

$p=139$, $q=2685619$,  $\zeta(x)=x^3+x+7$, $\tilde{c}=[436598,2118]$

$p=149$, $q=3307949$, $\zeta(x)=x^3+x+14$, $\tilde{c}=[1861398,3141536,1357853,1]$

$p=151$, $q=3442951$,  $\zeta(x)=x^3+x+5$, $\tilde{c}=[1492322,66]$

$p=163$, $q=4330747$,  $\zeta(x)=x^3+x+4$, $\tilde{c}=[2015256,4602]$

$p=167$, $q=4657463$,  $\zeta(x)=x^3+x+3$, $\tilde{c}=[3183263,109]$

$p=179$, $q=5735339$,  $\zeta(x)=x^3+x+4$, $\tilde{c}=[2740965,1219]$

$p=191$, $q=6967871$,  $\zeta(x)=x^3+x+3$, $\tilde{c}=[4789910,1160]$

$p=199$, $q=7880599$,  $\zeta(x)=x^3+x+13$, $\tilde{c}=[3457494,2368]$

$p=211$, $q=9393931$,  $\zeta(x)=x^3+x+24$, $\tilde{c}=[5457264,1168]$

$p=223$, $q=11089567$,  $\zeta(x)=x^3+x+9$, $\tilde{c}=[4722613,4305]$

$p=227$, $q=11697083$,  $\zeta(x)=x^3+x+9$, $\tilde{c}=[9051956,1442]$

$p=239$, $q=13651919$,  $\zeta(x)=x^3+x+11$, $\tilde{c}=[1597504,5918]$

$p=251$, $q=15813251$,  $\zeta(x)=x^3+x+7$, $\tilde{c}=[9285089,11965]$

$p=263$, $q=18191447$,  $\zeta(x)=x^3+x+8$, $\tilde{c}=[8313030,2840]$

$p=271$, $q=19902511$,  $\zeta(x)=x^3+x+4$, $\tilde{c}=[6563520,170]$

$p=283$, $q=22665187$,  $\zeta(x)=x^3+x+24$, $\tilde{c}=[2245440,3574]$

\bigskip

The following produce P1Fs for complete graphs of orders $q+1=p^5+1$, where $p$
is prime. They can be considered as extra entries for Table 7 in \cite{cyclatom}.

$p=19$, $q=2476099$,  $\zeta(x)=x^5+x+9$, $\tilde{c}=[949007,791]$

$p=23$, $q=6436343$,  $\zeta(x)=x^5+x+3$, $\tilde{c}=[1045440,7580]$

\subsection*{Acknowledgements}

This research was supported by the Monash eResearch Centre and
eSolutions-Research Support Services through the use of the MonARCH
HPC Cluster. The authors are also grateful to Mariusz Meszka and Alex Rosa for
helpful discussion, including pointers to the literature.


\begin{thebibliography}{99}

\bibitem{BMW06}
D.~Bryant, B.~Maenhaut and I.\,M.~Wanless, 
New families of atomic Latin squares and perfect one-factorisations,
{\it J.\ Combin.\ Theory Ser.\ A\/}, {\bf113} (2006), 608--624.

\bibitem{CC80}
  C.\,J.~Colbourn and M.\,J.~Colbourn,
  Combinatorial isomorphism problems involving 1-factorizations,
  \emph{Ars Combin.} {\bf9} (1980), 191--200.
  
\bibitem{DG96} J.\,H.~Dinitz and D.\,K.~Garnick, 
There are 23 non-isomorphic perfect one-factorisations of $K_{14}$, 
{\it J.~Combin.\ Des.\/} {\bf4} (1996), 1--4.

\bibitem{KMOW14}
P.~Kaski, A.\,D.\,S. Medeiros, P.\,R.\,J.~\"Osterg\aa rd and I.\,M.~Wanless, 
Switching in one-factorisations of complete graphs,
{\it Electron. J. Comb.\/} {\bf21}(2) (2014), \#P2.49.

\bibitem{Mes20}
M.~Meszka,
There are 3155 nonisomorphic perfect 1-factorizations of $K_{16}$,
submitted for publication.

\bibitem{MR03}
M.~Meszka and A.~Rosa, 
Perfect 1-factorizations of $K_{16}$ with nontrivial automorphism group, 
{\it J. Combin. Math. Combin. Comput.} {\bf47} (2003), 97--111.

\bibitem{Pik19}
D.\,A.~Pike, A perfect 1-factorisation of $K_{56}$,
\emph{J.\ Combin.\ Des.} {\bf27} (2019), 386--390.
  
\bibitem{Sea91}
E.~Seah, Perfect one-factorizations of the complete graph---a survey,
{\em Bull.\ Inst.\ Combin.\ Appl.\/} {\bf1} (1991), 59--70.

\bibitem{Wal97} W.\,D.~Wallis, {\it One-factorizations}, Kluwer Academic, 
Dordrecht, Netherlands, 1997.

\bibitem{cyclatom}
I.\,M.~Wanless, 
Atomic Latin squares based on cyclotomic orthomorphisms, 
{\it Electron.\ J.\ Combin.\/} {\bf12} (2005), R22.

\bibitem{catalog}
I.\,M.~Wanless, Author's homepage,
\url{http://users.monash.edu.au/~iwanless/data/P1F/newP1F.html}


\bibitem{syminherit}
I.\,M.~Wanless and E.\,C.~Ihrig, Symmetries that
Latin Squares Inherit from $1$-Factor\-izations,
{\it J.\ Combin.\ Designs\/} {\bf13} (2005), 157--172.

\bibitem{Wol09}
  A.\,J.~Wolfe,
  A perfect one-factorization of $K_{52}$,
  \emph{J. Combin. Des.} {\bf17} (2009), 190--196.

\end{thebibliography}
\end{document}